\documentclass[11pt,leqeq]{article}
\usepackage{amsmath,amsthm,amssymb,amsfonts,amscd, array}
\usepackage[all]{xy}

\newtheorem{thm}{Theorem}[section]
 \newtheorem{cor}[thm]{Corollary}
 
 \newtheorem{prop}[thm]{Proposition}
 \newtheorem{defn}[thm]{Definition}

\usepackage[dvips]{graphicx}
\DeclareFontFamily{OMS}{rsfs}{\skewchar\font'60}
\DeclareFontShape{OMS}{rsfs}{m}{n}{<-5>rsfs5 <5-7>rsfs7 <7->rsfs10
}{} \DeclareSymbolFont{rsfs}{OMS}{rsfs}{m}{n}
\DeclareSymbolFontAlphabet{\scr}{rsfs}

\def\Z{{\Bbb Z}}
\newcommand{\ve}{{\rm vect}}
\newcommand{\des}{\displaystyle}
\newcommand{\opa}{{\rm{OPar}}}
\date{}

\begin{document}
\title{Rota-Baxter Categories}
\author{Edmundo Castillo and Rafael D\'iaz}
\maketitle

\begin{abstract}
We introduce Rota-Baxter categories and construct examples of such
structures.\\


\end{abstract}

This work takes part in the efforts to understand the
categorification of rings and other related algebraic structures.
The idea of categorification of algebraic structures has been
around for several decades and has gradually become better
appreciated and understood. The expanding scope and applications
of the notion of categorification has been greatly influenced by
the works of Baez-Dolan
\cite{BaezDolan, BaezDolan3}, Crane-Frenkel
\cite{crane}, Crane-Yetter \cite{cy}, Khovanov
\cite{kho}, among others. The basic idea is that it is worthwhile to look
at the categorical foundations of set theoretical structures.
Often sets arise as the equivalences classes of  objects in a
category. Going from a category to the set of equivalences classes
of its objects is the process of decategorification.
Categorification goes in the reverse direction, uncovering
categories whose set of equivalences classes of objects reproduces
a given set. Categorifications always exist but are no unique.
Thus two general problems arise: the classification of
categorifications and the extraction of information regarding a
given set theoretical construction from its
categorical counterpart. \\

Our approach to the categorification of rings, reviewed in Section
\ref{caterings}, was first
discussed in
\cite{RDEP} with a view towards the categorification of the
ring of functions on non-commutative spaces and the
categorification of the algebra of annihilation and creation
operators. Further developments aimed at the elaboration of a
general setting for the study of the combinatorial properties of
rational numbers were reported in
\cite{BD1, BD2}, where the combinatorics of Bernoulli
numbers and hypergeometric functions, respectively, are discussed.
Our aim in this work is to study the categorification of
Rota-Baxter rings \cite{GCRota}, an algebraic structure under
current active research because of its capacity to unify notions
coming from probability theory, combinatorics, symmetric
functions, and the renormalization of Feynman integrals among
others. Applications of Rota-Baxter categories in the context of
renormalization of Feynman integrals will be study in the
forthcoming works \cite{cast, DP1}.\\

This paper is organized as follows. In Section \ref{caterings} we
introduce Rota-Baxter categories and provide a couple of basic
examples. In Section \ref{erbc} we provided the simplest and most
ubiquitous example of Rota-Baxter category. This sort of
Rota-Baxter category are constructed from any distributive
category and may be thought as categorifications of the
Rota-Baxter ring of formal Laurent series. In Section \ref{cmc} we
construct Rota-Baxter categories associated with an arbitrary
comonoidal category and a given Rota-Baxter category. This
construction should be thought as the categorification of the ring
of formal power series with coefficients in a Rota-Baxter ring. In
Sections \ref{ibf} and \ref{bf} we construct Rota-Baxter
categories from idempotent and arbitrary  bimonoidal functors,
respectively. In Sections
\ref{ci} and \ref{hoy} we define categorical integration and show
in three different contexts, categorical Riemannian integration,
discrete analogues of integration and categorical Jackson
integrals from $q$-calculus, that functorial integration provides
examples of (twisted) Rota-Baxter categories of various weights.
In Section \ref{cft} and \ref{qft} we construct Rota-Baxter
categories naturally arising from classical and quantum field
theory, respectively.

\section{Categorification of rings}\label{caterings}
We assume the reader to be familiar with basic notions of category
theory \cite{SMacLane}. Let us begin recalling the notion of
categorification of rings and semi-rings from
\cite{BD2, RDEP}.

\begin{defn}{\em A category $C$ is distributive if
it is equipped with functors $\oplus\colon C \times C
\rightarrow C$ and $\otimes\colon C \times C
\rightarrow C$ called sum and product, respectively.
There are distinguished objects $0$ and $1$ in $C$; $(C,\oplus,0)$
is a symmetric monoidal category with unit $0$; $(C, \otimes,1)$
is a  monoidal category with unit $1$. There are natural
isomorphisms $$x\otimes(y\oplus z)\simeq (x
\otimes y) \oplus (x \otimes z) \mbox{\ \  and \ \ } (x\oplus y)\otimes z \simeq
(x \otimes z) \oplus (y \otimes z),$$ for $x,y,z$ objects of $C$.
A distributive category have negative objects if it comes with a
functor $-\colon C \rightarrow C$ and for $x,y$ objects of $C$
there are natural isomorphisms
$$-(x\oplus y)\simeq -x\oplus -y,\ \  -0\simeq 0, \mbox{\ \ and \ \ } -(-x)\simeq
x.$$ }
\end{defn}

 Coherence theorems for distributive
categories were studied by Laplaza \cite{l1}. An interesting
research problem is to find coherence theorems for distributive
categories with negative objects.

\begin{defn}{\em
Let $C$ be a distributive category. A functor $P:C \rightarrow C$
is additive if for $x,y$ objects of $C$ there are natural
isomorphisms $P(x \oplus y)
\simeq P(x) \oplus P(y)$. If $C$ has a negative functor we also demand the existence of
natural isomorphisms $P(-x)\simeq -P(x)$. $P$ is bimonoidal if it
is additive and in addition there are natural isomorphisms
$P(x\otimes y)\simeq P(x) \otimes P(y)$.}
\end{defn}

\begin{defn} {\em A categorification  of a ring $R$ is a
distributive category $C$ with negative functor together with a
valuation map $|\ \ |\colon Ob(C) \rightarrow R$ such that:
$$|x|=|y| \mbox{\ \  if \ \ } x
\simeq y,\ \  |x\oplus y|=|x|+|y|,\ \  |x\otimes y| = |x||y|,\ \
|1|=1, \ \ |0|=0, \mbox{\ \ and \ \ } |-x|=-|x|.$$}
\end{defn}

If we omit the existence of the negative functor in the definition
above we arrive to the notion of categorification of semi-rings,
which will be used quite often in this work. Next we introduce the
main concept of this work, the notion of Rota-Baxter categories. A
Rota-Baxter ring, see
\cite{GCRota} and the references
therein, is a triple $(R,\lambda,p)$ where $R$ is a ring, $\lambda
\in
\{-1,0,1\}$, and $p\colon R\longrightarrow R$ is a morphism of
abelian groups satisfying: $$p(x)p(y)=p(xp(y))+p(p(x)y)+\lambda
p(xy).$$  $R$ may or may not have a unit, and may or may not be
commutative. Notice that the notion of Rota-Baxter semi-ring makes
perfect sense; for $\lambda=-1$, the required identity is
$$p(x)p(y)+ p(xy)=p(xp(y))+p(p(x)y).$$

\begin{defn}
{\em A Rota-Baxter category  of  weight $\lambda \in \{-1,0,1 \}$
is a distributive category $C$ together with an additive functor
$P: C \longrightarrow C$ and natural isomorphisms
\begin{eqnarray*}
P(x)\otimes P(y) \oplus P(x\otimes y)  &\simeq& P(P(x)\otimes y)
\oplus P(x\otimes P(y))
\\
P(x) \otimes P(y)&\simeq& P(P(x)\otimes y)
\oplus P(x\otimes P(y))\\
P(x) \otimes P(y)&\simeq& P(P(x)\otimes y)
\oplus P(x\otimes P(y))
\oplus  P(x\otimes y).
\end{eqnarray*}
for $x,y$ objects of $C$ and $\lambda=-1,0,1,$ respectively.}
\end{defn}

\begin{defn}{\em
A categorification of a Rota-Baxter ring $(R,p)$ is a Rota-Baxter
category $(C,P)$ together with a valuation map $|\ \ |\colon Ob(C)
\rightarrow R$ such that $|P(x)|= p(|x|)$ for $x$  object of
$C$.}
\end{defn}

A Rota-Baxter ring may be regarded as a Rota-Baxter category and,
as such, it is a categorification of itself. In the next sections
the reader will find interesting examples of Rota-Baxter
categories making this notion worth studying; we begin pointing
out a couple of simple but useful examples. Any ring may be
regarded as a Rota-Baxter ring with vanishing $p$. Any abelian
group $R$ provided with a group morphisms $P:R
\rightarrow R$ may be regarded as a Rota-Baxter ring with
multiplication constantly equal to zero. The analogues of these
simple facts hold in the categorical context as well.

\begin{prop}{\em  A distributive category may be regarded as a Rota-Baxter
category with functor $P$ constantly equal to zero. A distributive
category may be regarded as a Rota-Baxter category of weight $-1$
with $P$ equal to the identity functor. A symmetric monoidal
category $C$ together with an additive functor $P:C \rightarrow C$
may be regarded as a Rota-Baxter category with $\otimes$
constantly equal to zero.}
\end{prop}

Rota-Baxter rings with vanishing $P$ play a fundamental role in
the theory of renormalization as formalized by Connes and Kreimer
\cite{ck}. Rota-Baxter rings with vanishing product, though less
studied, should not be overlooked. If $(C,P)$ is a Rota-Baxter
category then we let $Ker(P)$ be the full subcategory of $C$ such
that $c$ is an object of $Ker(P)$ iff $P(c) \simeq 0$. Similarly,
let $Im(P)$ be the full subcategory of $C$ whose objects are
isomorphic to objects of the form $P(c)$ for some $c
\in Ob(C)$. The axioms for Rota-Baxter categories imply the
following result.
\begin{prop}{\em
$(Ker(P),0)$ is a Rota-Baxter category. $(Im(P),I)$ is a
Rota-Baxter category.}
\end{prop}

\section{Main examples}\label{erbc}

\begin{figure}[ht]
\begin{center}
\includegraphics[height=3cm]{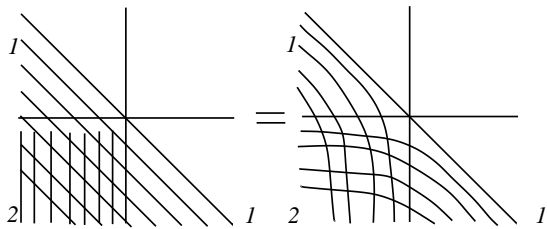}
\caption{ \ Geometric meaning of the Rota-Baxter identity. \label{bra1}}
\end{center}
\end{figure}

The reason why the examples considered in this section are
Rota-Baxter categories is succinctly  encoded in the identity
between sets with multiplicities shown in Figure \ref{bra1}.
Let us start with the simplest and most prominent example. Let
$\Z$-$\ve_{b}$ be the category of $\Z$-graded vector spaces
$V=\bigoplus V_n$ such that: $V_n$ is finite dimensional and there
exists $k \leq 0$ such that $V_{n}=0$ for $n\leq k.$
$\Z$-$\ve_{b}$ is a distributive category with direct sums and
tensor products given as usual by $$(V \oplus W)_n= (V_{n}\oplus
W_{n}) \mbox{ \ \ and \ \ } (V \otimes W)_n=
\bigoplus_{k+l=n}(V_{k}\otimes W_{l}).$$ Let $\mathbb{N}[z^{-1},z]]$ be the semi-ring of formal Laurent
series with integral coefficients, it is known that the map
$p:\mathbb{N}[z^{-1},z]] \rightarrow
\mathbb{N}[z^{-1}]$ turns $\mathbb{N}[z^{-1},z]]$ into a
Rota-Baxter semi-ring of weight $-1$.

\begin{prop}\label{vs}
{\em The functor $P\colon
\Z$-$\ve_{b}\rightarrow \Z$-$\ve_{b}$ given for $k<0$ by
$$P\Big(\bigoplus_{k \leq n}V_{n}\Big)=
\bigoplus_{k \leq n<0}V_{n},$$ turns
$\Z$-$\ve_{b}$ into a Rota-Baxter category of weight $-1$.
$\Z$-$\ve_{b}$ is a categorification of $\mathbb{N}[z^{-1},z]]$
with valuation map $| \ \ |:\Z$-$\ve_{b}\rightarrow
\mathbb{N}[z^{-1},z]]$ given by $|V|=\sum_n dim(V_n)z^n$.}
\end{prop}

Proposition \ref{vs} is an instance of a general construction of
Rota-Baxter categories to be developed presently. For a
distributive category $C$ let $C_b^{\Z}$ be the category whose
objects are maps $f: \Z
\rightarrow C,$ such that there exists $k \leq 0$ with $f(n)\simeq 0$ for $n<k$.
Morphisms in $C_b^{\Z}$ are given by $$C_b^{\Z}(f,g)=\prod_{n \in
\Z}C(f(n),g(n)).$$ The category $C_b^{\Z}$ is distributive category with
$\oplus$, $\otimes$, and negative functor given by: $$(f\oplus
g)(n)=f(n)\oplus g(n), \ \ (f\otimes
g)(n)=\bigoplus_{k+l=n}f(k)\otimes g(l), \mbox{\ \ and \ \
}(-f)(n)=-f(n),$$ where $f,g$ belong to $C_b^{\Z}$ and $k,l,n
\in
\Z.$
\begin{thm}\label{vt}
{\em The category $C_b^{\Z}$ together with the functor $P:C_b^{\Z}
\rightarrow C_b^{\Z}$ given by $P(f)=f_{<0}$ where
$$f_{<0}(n)=\left\{\begin{array}{cc}
0 & \mathrm{if}\ n \geq 0\\
f(n) & \mathrm{if}\ n <0
\end{array}\right.$$
is a Rota-Baxter category of weight $-1$. If $C$ is a
categorification of $R$, then $C_b^{\Z}$ is a categorification of
$R[z^{-1},z]].$}
\end{thm}

\begin{proof}
For $f$ and $g$ objects of $C_b^{\Z}$, we have that:
\begin{eqnarray*}
P( f\otimes g)&=& (f\otimes g)_{<0}\\
P(f)\otimes P(g)&=& f_{<0} \otimes g_{<0}\\
P(f\otimes P(g))&=& (f\otimes g_{<0})_{<0}\\
P( P(f)\otimes g) &=& (f_{<0}\otimes g)_{<0}.
\end{eqnarray*}
We have to check that there are canonical isomorphisms
$$(f\otimes g)_{<0} \oplus f_{<0} \otimes g_{<0} \simeq
(f_{<0}\otimes g)_{<0} \oplus  (f\otimes g_{<0})_{<0},$$ which we
do evaluating both sides at $n \in \Z$. If $n \geq 0$ we obtain
the identity $0 \oplus 0 = 0 \oplus 0$. If $n < 0$, then we have
to show that there are canonical isomorphisms
$$\bigoplus_{k+l=n}f(k)\otimes g(l) \ \ \bigoplus \bigoplus_{k+l=n, k <0, l<0}f(k)\otimes g(l) \simeq
\bigoplus_{k+l=n, k <0}f(k)\otimes g(l) \ \ \bigoplus \bigoplus_{k+l=n, l<0}f(k)\otimes g(l),$$
which is clear. This proves the first statement of the theorem.
For the second statement consider  the valuation map $| \ \ |:
C_b^{\Z}
\rightarrow  R[z^{-1},z]]$  given  by $$|f|=\sum_{n \in \Z}|f(n)|z^n$$
satisfies all required axioms.
\end{proof}

Proposition \ref{vs} is obtained from Theorem \ref{vt} letting $C$
be $vect$ the category of finite dimensional vector spaces with
valuation map $|V|=dim(V).$\\

Let $\mathbb{Z}$-$set$ be the category of $\mathbb{Z}$-graded
finite sets, i.e. an object of $\mathbb{Z}$-$set$ is a pair
$(x,f)$ where $x$ is a finite set and $f:x\rightarrow \mathbb{Z}$
is a map. Morphisms in $\mathbb{Z}$-$set$ from $(x,f)$ to $(y,g)$
are maps $\alpha:x
\rightarrow y$ such that $g \circ \alpha = f.$ Disjoint union and Cartesian product are given, respectively, by
$$(x,f)\sqcup (y,g)=(x\sqcup y,f \sqcup g) \mbox{ and } (x,f)\times
(y,g)=(x\times y,f\circ \pi_{x}+g\circ
\pi_{y}),$$ where $\pi_{x}$ and $\pi_{y}$ are the canonical
projections of $x \times y$ onto $x$ and $y$, respectively.
Consider the functor $P:\mathbb{Z}$-$set \longrightarrow
\mathbb{Z}$-$set$ given by
$$P(x,f)=(f^{-1}(-\infty,0),f|_{f^{-1}(-\infty,0)}).$$

\begin{prop}{\em
$(\mathbb{Z}$-$set,P)$ is a Rota-Baxter category of weight $-1$.
$\mathbb{Z}$-$set$ is a categorification of
$\mathbb{N}[z^{-1},z].$}
\end{prop}
\begin{proof}
The result follows from Theorem \ref{vt}  since there is a natural
functor $i:\mathbb{Z}$-$set \rightarrow set_b^{\Z}$ which exhibits
$\mathbb{Z}$-$set$ as a full subcategory of $set_b^{\Z}$ closed
under sum, product and $P$. The valuation map on $set_b^{\Z},$
induced from the valuation map on $set$ sending $x$ into its
cardinality $|x|$, restricts to a valuation map on
$\mathbb{Z}$-$set$.
\end{proof}

\section{Comonoidal categories}\label{cmc}

If $C$ is a co-ring, i.e., an abelian group provided with a
co-product, and $R$ is a Rota-Baxter ring, then the set $Hom(C,R)$
of morphisms of abelian groups from $C$ to $R$ is Rota-Baxter ring
with product $fg(x)=(f\otimes g)\Delta(x)$ and operator $p$ given
by $p(f)(c)=p(f(c))$ for $f
\in Hom(C,R)$ and $c \in C.$ We proceed to state the corresponding
facts for Rota-Baxter categories.\\

Let $Cat$ be the category whose objects are essentially small
categories, morphisms in $Cat$ are functors. We a define a functor
$\otimes: Cat \times Cat \longrightarrow Cat$ as follows. The
tensor product category $C \otimes D$  of the categories $C$ and
$D$ has as objects triples $(x,f,g)$ where $x$ is a finite set,
$f:x
\rightarrow Ob(C)$ and $g:x \rightarrow Ob(D)$ are
maps. Morphisms from $(x_1,f_1,g_1)$  to $(x_2,f_2,g_2)$, objects
of $C \otimes D$, are given by

$$C \otimes D((x_1,f_1,g_1),(x_2,f_2,g_2))=
\bigsqcup_{\alpha:x_1 \rightarrow x_2 }\prod_{i \in x}C(f_1(i),f_2(\alpha(i)))\times D(g_1(i),g_2(\alpha(i))),$$
where $\alpha:x_1 \rightarrow x_2$ is an arbitrary bijection.\\

The following definition formalizes the categorical analogue of
the notion of a co-ring without co-unit.

\begin{defn} {\em A category $D$ is comonoidal if it comes equipped with a
functor $\delta: D \longrightarrow D \otimes D$ such that there is
a natural isomorphisms $(\delta \otimes 1_D)\delta
\longrightarrow (1_D \otimes \delta)\delta$ satisfying McLane's pentagon axiom.}
\end{defn}

\begin{thm}\label{mas}
{\em Suppose $D$ with a functor $\delta: D \longrightarrow D
\otimes D$ is a comonoidal category and $C$ a Rota-Baxter
category. Then $C^D,$ the category of functors from $D$ to $C$, is
a Rota-Baxter category with functor $P$ given by
$P(F)(x)=P(F(x))$.}
\end{thm}

\begin{proof}
First we show that $C^D$ is distributive. Define sum and product
by $$(F+G)(x)=F(x)\oplus G(x)$$
$$(FG)(x)=\sum_{d\delta(x)}F(\delta_1(x)) \otimes
G(\delta_2(x)),$$ where $\delta(x): d\delta(x) \longrightarrow
Ob(D)\otimes Ob(D)$, where  $d\delta(x)$ is the domain of
$\delta(x)$. The negative functor is $(-F)(x)=-F(x).$ Next assume
$C$ is a Rota-Baxter category of weight $-1$, the other cases
being similar. The desired result follows from the natural
isomorphisms
$$P(F)\otimes P(G)(x) \oplus P(F\otimes G)(x) \simeq
\bigoplus_{d\delta(x)}P(F(\delta_1(x))
\otimes P(G(\delta_2(x))) \oplus P(F(\delta_1(x))
\otimes G(\delta_2(x))),$$
$$P(P(F)\otimes G)(x) \oplus P(F\otimes P(G)) \simeq
 \bigoplus_{d\delta(x)}P(F(\delta_1(x)))
\otimes G(\delta_2(x))  \oplus F(\delta_1(x)) \otimes P(G(\delta_2(x))).$$
\end{proof}

Given a positive integer $n$ we use the notation $[n]=\{1,2,\dots ,n\}.$ A category $D$ may be regarded as a comonoidal category with
functor $\delta: D \longrightarrow D \otimes D$ sending an object
$x$ in $D$ to the map $\delta(x):[1] \rightarrow D \times D$ such
that $\delta(x)(1)=(x,x).$ This canonical comonoidal structure
induces the monoidal structure on $C^D$ given by $FG(x)=F(x)G(x).$

\begin{cor}{\em
Asumme $C$ is a Rota-Baxter category. Then $C^D$ is a Rota-Baxter
category with functor $P:C^D
\rightarrow C^D$ given by $P(F)(x)=P(F(x))$, and the product of functors given by
$FG(x)=F(x)G(x).$}
\end{cor}

Let us consider a rather simple example of  comonoidal category.
Recall that a set $x$ may be regarded as the category with objects
$x$ and identity morphisms only. The category $[n]
\times [n]$ is comonoidal with functor $\delta: [n] \times
[n] \longrightarrow ([n] \times [n]) \otimes ([n]
\times [n])$  such that $\delta(i,j)$ is the map $$\delta(i,j):[n]
\longrightarrow ([n] \times [n]) \otimes ([n]
\times [n])$$ given by $$\delta(i,j)(k)=((i,k), (k,j)). $$
If $C$ is a distributive category then  $M_{n}(C)=C^{[n]\times
[n]}$, the category of $n \times n$ matrices with values in $C$,
is also a distributive category. Concretely, an object $A$ in
$M_{n}(C)$ is a family $A_{ij}$ of objects in $C$, for $1\leq i,j
\leq n$. Sum and product of objects in  $M_{n}(C)$ are given, respectively, by $$(A
\oplus B)_{ij}=A_{ij}\oplus B_{ij} \mbox{\ \  and \ \ }(A\otimes
B)_{ij}=\des \bigoplus_{k}(A_{ik}\otimes B_{kj}),$$ for $A,B \in
M_{n}(C).$

\begin{cor}{\em If $C$ is a distributive category then $M_{n}(C)$ is a
distributive category. If $C$ is a Rota-Baxter category, then
$M_{n}(C)$ is a Rota-Baxter category. If $C$ is a categorification
of $R$, then $M_{n}(C)$ is a categorification of $M_{n}(R)$.}
\end{cor}

For our next constructions we need the theory of species
introduced by Joyal
\cite{j2} and further elaborated by Bergeron, Labelle, and Leroux
\cite{Bergeron}. Let $\mathbb{B}^n$ be the category whose objects
are pairs $(x,f)$ where $x$ is a finite set and $f:x \rightarrow
[n]$ is any map. Morphisms in $\mathbb{B}^n$ from $(x,f)$ to
$(y,g)$ are maps $\alpha:x \rightarrow y$ such that $g \circ
\alpha =f$. For a distributive category $C$ we let
$C^{\mathbb{B}^n},$ the category of $C$-species in $n$ variables,
be the category of functors from $\mathbb{B}^n$ to $C$.\\

The category $\mathbb{B}^n$ is comonoidal  with functor $\delta:
\mathbb{B}^n
\rightarrow \mathbb{B}^n \otimes \mathbb{B}^n$ given on  $(x,f)$ in $\mathbb{B}^n$ by
the map $$\delta(x,f):Par_2(x) \rightarrow Ob(\mathbb{B}^n)
\otimes Ob(\mathbb{B}^n)$$ such that
$$\delta(x,f)(x_1,x_2)=((x_1,f|_{x_{1}}),(x_2, f|_{x_2})),$$ where
$Par_2(x)$ is the set of pairs $(x_1,x_2)$ such that $x_1
\sqcup x_2 = x$. It follows that $C^{\mathbb{B}^n}$ is a distributive
category with sum and product given by $$(F+G)(x,f)=F(x,f)\sqcup
G(x,f)\mbox{\ \ and \ \ }(FG)(x,f)=\bigoplus_{x_1\sqcup x_2
=x}F(x_1, f|_{x_1})
\otimes G(x_2,f|_{x_2}).$$
If $R$ is a ring then we let $R[[x_1,...,x_n]]$ be  the ring of
formal divided power series in variables $x_1,...,x_n$. The latter
algebra is the free $R$-module generated by symbols:
$$\frac{x^k}{k!}\mbox{\ \  where \ \ } k
\in \mathbb{N}^n,\ \  x^k=x_1^{k_1}...x_n^{k_n},  \mbox{\ \ and \ \ }
k!=k_1!...k_n!.$$  The product is defined on generators via the
identity
$$\frac{x^k}{k!}\frac{x^s}{s!}=
\binom{k+s}{s} \frac{x^{k+s}}{(k+s)!},$$ where
$\binom{k+s}{s}=\prod_{i=1}^n\binom{k_i + s_i}{s_i}$.

\begin{cor}
{\em  If $(C,P)$ is Rota-Baxter category $\lambda$, then
$C^{\mathbb{B}^n}$ is a Rota-Baxter category of weight $\lambda$
with functor $P$ given by $P(F)(x,f)= P(F(x,f)),$ for $F$ in
$C^{\mathbb{B}^n}$ and $(x,f)$ in $\mathbb{B}^n$. If $C$ is a
categorification of $R$, then $C^{\mathbb{B}^n}$ is a
categorification of $R[[x_1,...,x_n]]$.}
\end{cor}
\begin{proof}
Follows from Theorem \ref{mas}. The valuation map $|\
\ |:C^{\mathbb{B}^n}
\rightarrow R[[x_1,...,x_n]]$ is defined by $$F =
\sum_{k \in \mathbb{N}^n}F([k])\frac{x^k}{k!}\mbox{\ \  where \ \ } [k]=([k_1],...,[k_n]).$$
\end{proof}

Next we consider non-commutative species introduced in
\cite{RDEP}. Let $R\langle
\langle x_1,\dots, x_n
\rangle
\rangle$ be the ring of formal power series in non-commutative
variables $x_1,\dots,x_n$ and coefficients in $R$. We construct a
categorification $C^{L_n}$ of $R \langle \langle x_1,\dots, x_n
\rangle \rangle$ with the property that each valuation $|\ \ |\colon
Ob(C)\rightarrow R$, induces a valuation $$|\ \ |\colon Ob(
C^{L_n})
\rightarrow  R\langle
\langle x_1,\cdots, x_n \rangle \rangle .$$  Let $L_n$ be the category
whose objects are triples $(x,<,f)$ where $x$ is a finite set, $<$
is a linear order on $x$, $f\colon x
\rightarrow [n]$ is a map. Morphisms  from
$(x,<,f)$  to $(y,<,g)$ are given by
$$L_m((x,<,f), (y,<,g))= \{   \varphi\colon x
\rightarrow y \ |  \ g\circ \varphi=f,\ \mbox{and}\ \varphi(i)<
\varphi(j)\ \ \mbox{for all}\ \ i<j \}.$$ The disjoint union
$(x_1,<_1)\sqcup (x_2,<_2)$ of linearly order sets is  $(x_1\sqcup
x_2, <)$, where the order on $x_1\sqcup x_2$ extends the order on
$x_1$, the order on $x_2$, and $i < j$ for $i\in x_1$, $j\in x_2$.
An order partition in $n$-blocks of $(x,<)$ is a $n$-tuple
$(x_1,<_1),\cdots,(x_n,<_n)$ of posets such that $(x_1,<)\sqcup
\cdots \sqcup (x_n,<)=(x,<)$. Let $OPar_n(x,<)$ be
the set of order partitions of $(x,<)$ in $n$ blocks. $L_n$ is
comonoidal category with $\delta:L_n \rightarrow L_n
\otimes L_n$ sending $(x,<,f)$ into  the map $$\delta(x,<):
OPar_2(x,<)\longrightarrow Ob(L_n) \times Ob(L_n)$$ such that
$$\delta(x,<)((x_1,<_1),(x_1,<_2))=(x_1,<_1,f{|}_{x_1}),(x_1,<_2,f{|}_{x_2}).$$
It follows that $C^{L_n}$ is distributive with sum $(F +
G)(x,<,f)=F(x,<,f)\oplus G(x,<,f)$ and product
$$(F G)(x,<,f)= \bigoplus F(x_1,<_1,f{|}_{x_1})\otimes
G(x_2,<_2,f{|}_{x_2}),$$ where the sum runs over all pairs
$((x_1,<),(x_2,<_1)) \in OPar_2(x,<_2).$

\begin{cor}{\em If $C$ is a categorification of a Rota-Baxter ring $R$, then
$C^{L_n}$  is a categorification of the Rota-Baxter ring $R
\langle
\langle x_1,\dots, x_n \rangle \rangle.$}
\end{cor}
\begin{proof} By Theorem \ref{mas} we only need to define the valuation map $| \ \ |:C^{L_n}
\rightarrow R \langle \langle x_1,\dots,
x_n \rangle \rangle$ which is given by $|F|=\sum_{f:[m]
\rightarrow [n]} |F([m],<,f)|x^f,$ where $x_f =
x_{f(1)}...x_{f(m)}$.

\end{proof}

\section{Idempotent bimonoidal functors}\label{ibf}

In this section we provided a general constructions which
generates a wide variety of examples of Rota-Baxter categories.

\begin{thm}\label{bo}{\em
Let $C$ be a distributive category and $P:C \rightarrow C$ a
bimonoidal functor such that there is a natural isomorphisms $P^2
\simeq P$. Then $(C,P)$ is a Rota-Baxter category of weight $-1$.}
\end{thm}

\begin{proof}
Since $P$ is bimonoidal and $P^2=P$ we have natural isomorphisms
$$P(P(x)\otimes y)\simeq P(x\otimes P(y))\simeq P(x\otimes y)
\simeq  P(x)\otimes P(y).$$
\end{proof}

\begin{thm}{\em Let $C_1$ and $C_2$ be distributive categories. $C_1 \times C_2$ is a distributive
category with sums, products, and negative functor given by
$(c_1,c_2)
\oplus (d_1,d_2)= (c_1 \oplus c_2, d_1 \oplus d_2)$, $(c_1,c_2)
\otimes (d_1,d_2)= (c_1 \otimes c_2, d_1 \otimes d_2)$ and $-(c_1,c_2)=(-c_1, -c_2).$ The functor $\pi_1:C_1\times C_2
\rightarrow C_1 \times C_2$ sending $(c_1,c_2)$ into $(c_1,0)$ turns $C_1 \times C_2$ into a Rota-Baxter category of weight $-1$. }
\end{thm}
\begin{proof}
Follows from Theorem \ref{bo} and the fact that $\pi_1$ is
bimonoidal and $\pi_1^2=\pi_1$. Equivalently, one can check the
identities $$\begin{array}{ccccc}
P((c_1, c_2)\otimes (d_1, d_2)) & = & P(c_1, c_2) \otimes P(d_1, d_2) & = & P(P(c_1, c_2)
\otimes (d_1, d_2)) \\
 & = & P((c_1, c_2) \otimes P(d_1, d_2)) & =& (c_1 \otimes
d_1, 0).
             \end{array}
$$
\end{proof}

We now consider a list of examples of Rota-Baxter categories all
of which are corollaries of Theorem \ref{bo}. Examples will become
gradually closer to geometric and topological matters, and we will
see that fundamental constructions in topology, such as the
intersection homology for stratified manifolds, can be naturally
recast within our settings. Let $f$ be a morphism in a category
$C$. We denote by $d(f)$ and $c(f)$ the domain and codomain of
$f$, respectively. If $x$ is an object of $C$ we let $1_x$ be the
identity morphisms from $x$ to itself, thus $d(1_x)=c(1_x)=x.$ We
let $C^{(1)}$ be the category whose objects are morphisms in $C$.
A morphisms in $C^{(1)}$ from $f$ to $g$ is a pair of morphisms
$(\alpha_1,
\alpha_2)$ in $C$ such that $\alpha_2
\circ f = g \circ \alpha_1$. If $C$ is a distributive category,
then we have induced functors $\oplus: C^{(1)} \times C^{(1)}
\rightarrow C^{(1)}$ and $\otimes: C^{(1)} \times C^{(1)}
\rightarrow C^{(1)}$ turning $C^{(1)}$ into a distributive
category.

\begin{thm}\label{imp}
{\em
If $C$ is a distributive category, then $C^{(1)}$ is a Rota-Baxter
category of weight $-1$ with functor $P:C^{(1)} \rightarrow
C^{(1)}$ given by $P(f)=1_{d(f)}$.}
\end{thm}

\begin{proof}
Follows from the natural isomorphisms $$P(f
\otimes g)\simeq P(f) \otimes P(g) \simeq P(P(f) \otimes g) \simeq P(f \otimes
P(g)) \simeq 1_{d(f)}\otimes 1_{d(g)}.$$
\end{proof}

Recall that a morphisms $f$ in a category $C$ is injective or
monic if $g \circ f= h \circ f$ implies $f=g$ for all morphisms
$f,g$ in $C^{(1)}.$ Let $I(C^{(1)})$ be the full subcategory of
$C^{(1)}$ whose objects are injective morphisms. Assume that $C$
is a distributive category and that the functors $\oplus$ and
$\otimes$ on $C^{(1)}$ induce by restriction functors $\oplus$ and
$\otimes$ on $I(C^{(1)})$. The functor $P$ on $C^{(1)}$ induces by
restriction a functor $P$ on $I(C^{(1)})$.
\begin{cor}
{\em $(I(C^{(1)}), P)$ is a Rota-Baxter category of weight $-1$.}
\end{cor}

Consider the category $Top_{\subseteq}$ whose objects are pairs of
topological spaces $(X_1,X_2)$ with $X_1 \subseteq X_2$. A
morphism form $(X_1,X_2)$ to $(Y_1,Y_2)$ is a continuous map
$f:X_2 \rightarrow Y_2$ such that $f(X_1) \subseteq Y_1.$
Componentwise disjoint union and Cartesian product give
$Top_{\subseteq}$ the structure of a distributive category. The
functor $P:Top_{\subseteq} \rightarrow Top_{\subseteq}$ given by
$P(X_1,X_2)=(X_1,X_1)$ is bimonoidal and $P^2=P.$
\begin{cor}
{\em $(Top_{\subseteq},P)$ is a Rota-Baxter category of weight
$-1$.}
\end{cor}

Consider $Vect_{\subseteq}$ the category  whose objects are pairs
$(V,W)$ where $W$ is a vector space and $V$ a subspace of $W$.
Morphisms are linear transformations between the bigger spaces
that respect  the given subspaces. $Vect_{\subseteq}$ is a
distributive category with sum and product given by componentwise
direct sum and tensor product. The functor $P:Vect_{\subseteq}
\rightarrow Vect_{\subseteq}$ given by $P(V,W)=(V,V)$ is
bimonoidal and satisfies $P^2=P$.

\begin{cor}
{\em $(Vect_{\subseteq},P)$ is a Rota-Baxter category of weight
$-1$.}
\end{cor}

Let $vect_{\subseteq}$ be the full subcategory of
$Vect_{\subseteq}$ whose objects are pairs of finite dimensional
vector spaces. Let $P: vect_{\subseteq}
\rightarrow vect_{\subseteq}$ be the functor given by $P(V,W)=(0,W/V)$, where
$W/V$ is the quotient vector space and $0$ the vector space with
one element. Again it is easy to check that $P$ is bimonoidal and
that $P^2=P.$

\begin{cor}
{\em $(vect_{\subseteq},P)$ is a Rota-Baxter category of weight
$-1$.}
\end{cor}

Let $\mathbb{T}$ be the category whose objects are triples
$(x,y,f)$ where $x
\subseteq y$ are sets and $f:y \rightarrow y$ is a map. Morphisms
in $\mathbb{T}$ from $(x_1,y_1,f_1)$ to $(x_2,y_2,f_2)$ are maps
$\alpha:y \rightarrow y$ such that $\alpha(x_1) \subseteq x_2$ and
$\alpha \circ f_1= f_2 \circ
\alpha.$ Sum and product are given by
$$(x_1,y_1,f_1) \sqcup (x_2,y_2,f_2) = (x_1 \sqcup x_2, y_1 \sqcup y_2, f_1 \sqcup f_2),$$
$$(x_1,y_1,f_1) \times (x_2,y_2,f_2) = (x_1 \times x_2, y_1 \times y_2, f_1 \times f_2).$$
The functor $P:\mathbb{T} \rightarrow \mathbb{T}$ given by
$P(x,y,f)=(a,a,g),$ where $a= \{i \in x \ | \  f(i) \in x \}$ and
$g= f|_a,$ is bimonoidal and $P^2=P.$
\begin{cor}{\em
$(\mathbb{T},P)$ is a Rota-Baxter category of weight $-1$.}
\end{cor}

An interesting refinement of $\mathbb{T}$ is the category
$\mathbb{T}_G$ where $G$ is a group. Objects in $\mathbb{T}_G$ are
triples $(x,y,\rho)$ where $x \subseteq y$ are sets and $\rho:G
\times x \rightarrow x$ is a group action of $G$ on $x$.
The distributive structure on $\mathbb{T}_G$ is defined just as
for $\mathbb{T}$. The functor $P:\mathbb{T}_G \rightarrow
\mathbb{T}_G$ given by $P(x,y,\rho)=(a,a,\rho|_a),$ where $a= \{i \in x \ | \
gi \in x \mbox{ for all } g \in G \}$ and $\rho|_a$ is the
restriction of the action of $\rho$ on $x$ to $a$, is bimonoidal
and $P^2=P.$

\begin{cor}{\em
$(\mathbb{T}_G,P)$ is a Rota-Baxter category of weight $-1$.}
\end{cor}

Consider the category $Comp_{\subseteq}$ whose objects are triples
$(V,W,d),$ where $V \subseteq W$ are  $\mathbb{Z}$-graded vector
spaces and $\partial:W \rightarrow W$ is a degree $-1$ linear map
such that $\partial \circ \partial =0$. Componentwise direct sum
and tensor product of $\mathbb{Z}$-graded vector spaces give
$Comp_{\subseteq}$  the structure of a distributive category. The
differential on the sum and tensor product is given by $d_{V\oplus
W} = d_V \oplus d_W$ and $d_{V\otimes W} = d_V
\otimes 1_W
\oplus 1_V \otimes d_W$. The homology of $(V,W,d)$ in
$Comp_{\subseteq}$ is by definition given by
$$H((V,W,d))=H(W,d)).$$ The functor $P:Comp_{\subseteq} \rightarrow
Comp_{\subseteq}$ given by $P(V,W,d)=(X,X,d_X),$ where for $i
\in
\mathbb{Z}$ we set $$X_i= \{v \in V_i \ | \ dv \in V_{i-1} \}$$ and
$d_X$ is the restriction of $d$ to $X$, is bimonoidal and $P^2=P$.

\begin{cor}\label{ih}{\em
$(Comp_{\subseteq},P)$ is a Rota-Baxter category of weight $-1$.}
\end{cor}

Let us point out the relation between the Corollary \ref{ih} above
and the notion of intersection homology as introduced by Goresky
and McPhearson in \cite{go}. In a nutshell the construction of
intersection homology may be summarized as follows. Fix a map
$p:\mathbb{Z}^{\geq 1} \rightarrow
\mathbb{Z},$ called the perversity, satisfying $$p(k) \leq p(k+1) \leq p(k)+1 \mbox{ \ \
and\ \ } p(1)=p(2)=0.$$ Construct functor $C^p:sman \rightarrow
Comp_{\subseteq}$. An objects in the groupoid $sman$ of stratified
manifold of dimension $n$ is a topological space $X$ together with
a filtration $$\emptyset = X_{-1}\subseteq X_0
\subseteq X_1 \subseteq ... \subseteq X_n=X$$ such that
$X_j \setminus X_{j-1}$ with the induced topology,
is a smooth manifold of dimension $j$. Morphisms are
homeomorphisms which are smooth when restricted to the smooth
pieces $X_j \setminus X_{j-1}$. \\

Given a stratified manifold and a perversity $p$, then $(C^p(X),
C(X))$ be the object of $Comp_{\subseteq}$ where $C(X)$ is the
complex of chains on $X$, i.e. $C_i(X)$ is the space of continuous
maps from $\Delta_i$ into $X$, and $C_i^p(X)$ is the subspace of
$C_i(X)$ generated by "allowed chains" of dimension $i$, i.e. the
space generated by chains $c:\Delta_i
\longrightarrow X$ such that  $c^{-1}(X_j
\setminus X_{j-1})$, for $j<n$, is included in the union of sub-simplices of
$\Delta_i$ of dimension less or equal to $i +j - n + p(j)$. The
intersection homology $IH^p(X)$ with perversity $p$ of a
stratified manifold $X$ is defined to be the homology of
$P(C^p(X))$, i.e. $IH^p(X)=H(P(C^p(X))).$ The discovery of the
intersection homology for stratified manifolds has been regarded
as one of the greatest achievements of twenty century mathematics.
One may only wonder at the fact that a Rota-Baxter functor was
already lurking around its very definition.

\section{Bimonoidal functors}\label{bf}

In this section we  assume that $C$ is a distributive category
with infinite sums and infinite distributivity.

\begin{thm}\label{loco}
{\em Let $F\colon C \to C$ be a bimonoidal functor. The functor
$P:C \rightarrow C$ given by
$$P(x)=\des\bigoplus_{n=0}^{\infty}F^{n}(x)$$ makes $C$ a
Rota-Baxter category of weight $-1$.}
\end{thm}
\begin{proof}
The desired result follows from the following natural isomorphisms
$$P(x)\otimes P(y) \simeq
\des\bigoplus_{n,m=0}^{\infty}(F^{n}(x)\otimes F^{m}(y)),\ \
P(x\otimes y)
\simeq  \des\bigoplus_{n=0}^{\infty}(F^{n}(x)\otimes F^{n}(y)),$$
$$P(x\otimes P(y)) \simeq
\des\bigoplus_{0=n\leq m}^{\infty}(F^{n}(x)\otimes F^{m}(y)), \ \
P(P(x)\otimes y)\simeq
\des\bigoplus_{0=m\leq n}^{\infty}(F^{n}(x)\otimes F^{m}(y)).$$

\end{proof}

\begin{cor}
{\em Fix species $F_i$ in $C_0^{\mathbb{B}^{n}}$ for $1 \leq i
\leq n$. The functor $(F_1,...,F_n): C_0^{\mathbb{B}^{n}}
\longrightarrow C_0^{\mathbb{B}^{n}}$ given by
$$(F_1,...,F_n)(F)=F\circ(F_1,...,F_n)$$ is bimonoidal. The functor
$P=P(F_1,...,F_n):C_0^{\mathbb{B}^{n}}
\longrightarrow C_0^{\mathbb{B}^{n}}$ given by $$P(F)=\sum_{m=0}^{\infty}F \circ (F_{1},\cdots ,F_{m})^{\circ m}$$
gives $C_0^{\mathbb{B}^{n}}$ the structure of a Rota-Baxter
category of weight $-1$.}
\end{cor}

Let $C_0^{L_n}$ be the full subcategory of $C^{L_n}$ whose objects
are functors $F$ such that $F(\emptyset)=0\in Ob(C)$.  Let $F,
F_1,\cdots ,F_n$ be non-commutative species in $C_0^{L_n}$ and
$(x,<,f)$ an object of $L_n$. The composition or substitution of
non-commutative species is given by
$$F(F_1,\cdots,F_n)(x,<,f)=\des{\bigoplus_{ p,g}F(p,<_{p},
g)\otimes \bigotimes_{B\in p} F_{p(B)}(B,<_{B}, f{|}_{B})},$$
where the sum runs over all $p\in\opa(x,<)$ and $g\colon
p\rightarrow [n]$.

\begin{cor}{\em
Fix species $F_i$ in $C_0^{L^{n}}$ for $1 \leq i
\leq n$. The functor $(F_1,...,F_n): C_0^{L^{n}}
\longrightarrow C_0^{L^{n}}$ given by
$$(F_1,...,F_n)(F)=F\circ(F_1,...,F_n)$$ is bimonoidal. The functor
$P=P(F_1,...,F_n):C_0^{L^{n}}
\longrightarrow C_0^{L^{n}}$ given by $$P(F)=\sum_{m=0}^{\infty}F \circ (F_{1},\cdots ,F_{m})^{\circ
m}$$
gives $C_0^{L^{n}}$ the structure of a Rota-Baxter category of
weight $-1$.}
\end{cor}

Let us close this section with an example related with
$q$-calculus. For a nice introduction to $q$-calculus the reader
may consult \cite{ck}. Recent results on $q$-calculus related with
Gaussian and Feynman integration are given in \cite{ED1, ED3,
CTT}. In Section \ref{hoy} we discuss further applications to
$q$-calculus. Consider the ring $R[[x,q]]$ of formal power series
in variables $x,q$ with coefficients in $R$. A fundamental role in
$q$-calculus is play by the shift operator $$s: R[[x,q]]
\longrightarrow R[[x,q]]$$ given by $$s(f)(x,q)=f(qx,q)$$ for $f \in
R[[x,q]].$ Suppose $C$ is a categorification of $R$, then
$C^{\mathbb{B}^2}$ is a categorification of $R[[x,q]]$. Our next
goal is to find a categorification of the shift operator, namely,
we define functor $S: C^{\mathbb{B}^2} \longrightarrow
C^{\mathbb{B}^2}$ such that $|S(F)|=s(|F|)$ for any $F$ in
$C^{\mathbb{B}^2}.$ Let $Inj:\mathbb{B}^2 \longrightarrow
\mathbb{B}$ be the species such that
$$Inj(x,j)=\{\alpha:x \rightarrow y \ | \ \alpha \mbox{ is
injective } \},$$ and define $S: C^{\mathbb{B}^2} \longrightarrow
C^{\mathbb{B}^2}$ by the rule $$S(F)(x,y)= \bigoplus_{\alpha \in
Inj(x,y)}F(x, y \setminus \alpha(x)).$$
\begin{thm}{\em The functor
$S: C^{\mathbb{B}^2} \longrightarrow C^{\mathbb{B}^2}$ given by
the formula above is bimonoidal and satisfies $|S(F)|=s(|F|)$. The
functor $P_S:C^{\mathbb{B}^2} \longrightarrow C^{\mathbb{B}^2}$
given by $$P_S(F)(x,y)=\bigoplus_{k\geq 0}
\bigoplus_{\alpha
\in Inj(x,y)^k}F(x, y \setminus \cup_{i=1}^k\alpha_i(x)),$$ where $\alpha=(\alpha_1,...\alpha_k)$ and
$\alpha_i(x) \cap \alpha_j(x)=\emptyset,$  turns
$C^{\mathbb{B}^2}$ into a Rota-Baxter category of
weight $-1$.}
\end{thm}
\begin{proof} The first part follows from the identities
$$|S(F)|= \sum_{k \leq
n}|S(F)([k],[n])|\frac{x^k}{k!}\frac{x^n}{n!}=
\sum_{k \leq
n}k! \binom{n}{k}f_{k,n-k}\frac{x^k}{k!}\frac{x^n}{n!}= s(|F|).$$
From Theorem \ref{loco} we know that the functor
$P_S(F)=\des\bigoplus_{n=0}^{\infty}S^n(F)$ is Rota-Baxter of
weight $-1$. It is easy to check that $P_S$ is given by the
formula above.
\end{proof}

From the expression above for $P_S$ we can easily compute
$$p_s(f)=|P_S(F)|=\sum_{k\leq n}(\sum_{pk\leq
n}\frac{n!F_{n,n-pk}}{(n-pk)!})\frac{x^k}{k!}\frac{x^n}{n!}.$$

\section{Functorial integration}\label{ci}

The most prominent example of a Rota-Baxter ring of weight $0$ is
the ring of continuous functions on the real line. The Rota-Baxter
operator is integration $$P(f)=\int_0^x fdt.$$ The Rota-Baxter identity in this case is equivalent
to the integration by parts formula. We
consider categorical analogues of this construction, what is
needed are categorifications of continuous functions such that it
is possible to define categorical analogues of the notion of
integration. Sinse we do not have at our disposal a
surjective categorification of the ring of continuous or smooth
functions on the real line, we restrict our attention to the
sub-ring of polynomial functions. Indeed, we work with
formal power series instead of polynomials functions.
Thus, we are looking for categorifications $| \ \
|:Ob(C)
\rightarrow R[[x]]$ of the ring of formal power series with
coefficients in a commutative ring $R$, with the property that
there exists a Rota-Baxter functor $P:C \rightarrow C$ such that
$|P(c)|=\int_0^x |c|dt.$ The categories with these properties  we
know of are categories of functors, and thus
$P$ provides a notion of functorial integration.\\

Let $L$ be the category of finite linearly order sets. Morphisms
in $L$ are order preserving bijections. If $x$ is a linearly order
set and $x_1
\sqcup x_2=x$, then $x_1$ and $x_2$ are linearly order with
the induced orders. If $x$ is a linearly order finite set, then
for $k \leq |x|$, we let $m_k(x)$ be the maximal interval of
length $k$ of $x$. For example $m(x)=m_1(x)$ is the maximal
element of $x$. Suppose that $C$ is a categorification of $R$ and
let $C^L$ be the category of functors from $L$ to $C$. Objects of
$C^L$ are called non-symmetric or linear $C$-species. Sum and
product of linear species are given by $F+G(x)=F(x)
\oplus G(x)$ and $$FG(x)=\bigoplus_{x_1 \sqcup x_2=x}F(x_1)\otimes
F(x_2),$$ for $F,G$
in $C^L$ and $x$ in $L$. Define functor $P:C^L
\rightarrow C^L$ by the rule $P(F)(x)=F(x\setminus m(x)).$

\begin{thm}{\em
$(C^L, P)$ is a Rota-Baxter category of weight $0$. The valuation
map $|\ \ |:C^L \rightarrow R[[x]]$ given by $$|F|=\sum_{n\in
\mathbb{N}}|F(n)|\frac{x^n}{n!}$$ satisfies $|P(F)|=\int_0^x |F|dt$ for any linear $C$-species $F$.}
\end{thm}

\begin{figure}[h]
\begin{minipage}[t]{0.5\linewidth}
\begin{center}
\includegraphics[height=2cm]{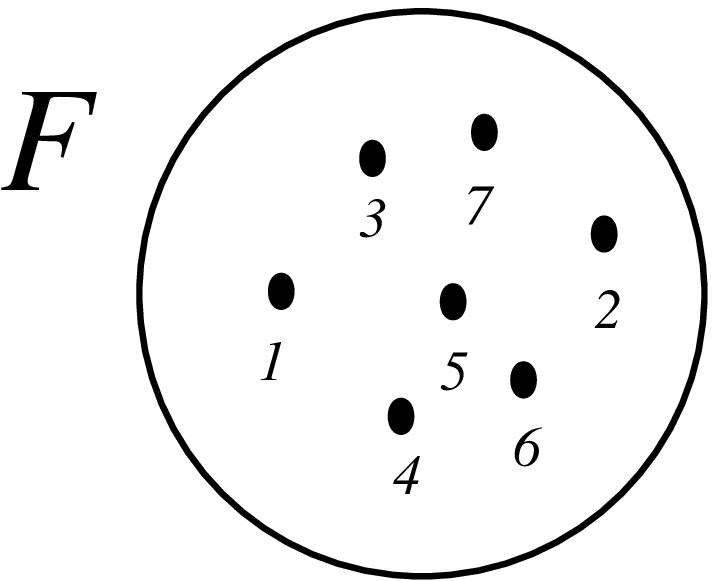}
\caption{ \ Evaluation of a species. \label{rb}}
\end{center}
\end{minipage}
\begin{minipage}[t]{0.5\linewidth}
\begin{center}
\includegraphics[height=2cm]{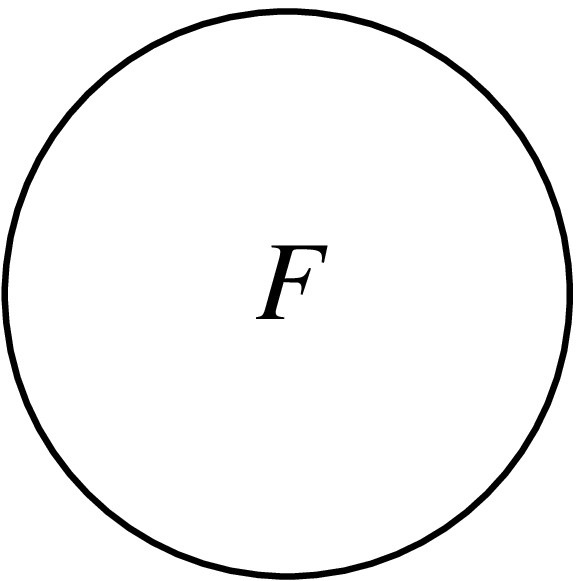}
\caption{ \ Abstract representation. \label{rb5}}
\end{center}
\end{minipage}
\end{figure}

There is a forgetful functor $f:L \rightarrow
\mathbb{B}$ which sends an ordered set $(x, <)$ into $x$. The functor $f$ induces by pullback the bimonoidal functor
$f^*: C^{\mathbb{B}} \rightarrow C^{L}$ and the commutative
diagram
\[\xymatrix @R=.3in  @C=.5in
{C^{\mathbb{B}}
\ar[r]^{| \  |} \ar[d]_{f^*} &
R[[x]] \ar[d]^{I}
\\ C^{L} \ar[r]^{| \  |} & R[[x]]} \]
Unfortunately the functor $P:C^{L} \rightarrow C^{L}$ is not well
defined on $C^{\mathbb{B}}$ since there is no canonical way to
choose an element from an unordered set. We see that in order to
define $P$ we should break the symmetries of $\mathbb{B}$, i.e.
reduce the isotropy groups to the identity and work with the
non-symmetric category $L$. The proof that $C^L$ is distributive
may be found in \cite{RDEP}. We show that $P$ is a Rota-Baxter
functor using graphical notation, the reader should bare in mind
that the grammatical codification of the pictorial proof is purely
mechanical. For example the evaluation of the species $F$ on the
set $x=\{1,2,3,4,5,6,7 \}$ is shown in Figure \ref{rb}. If we do
not need, and this is usually the case, to specify the elements of
$x$, then we prefer the use the abstract representation shown in
Figure \ref{rb5}.

\begin{figure}[h]
\begin{minipage}[t]{0.5\linewidth}
\begin{center}
\includegraphics[height=2cm]{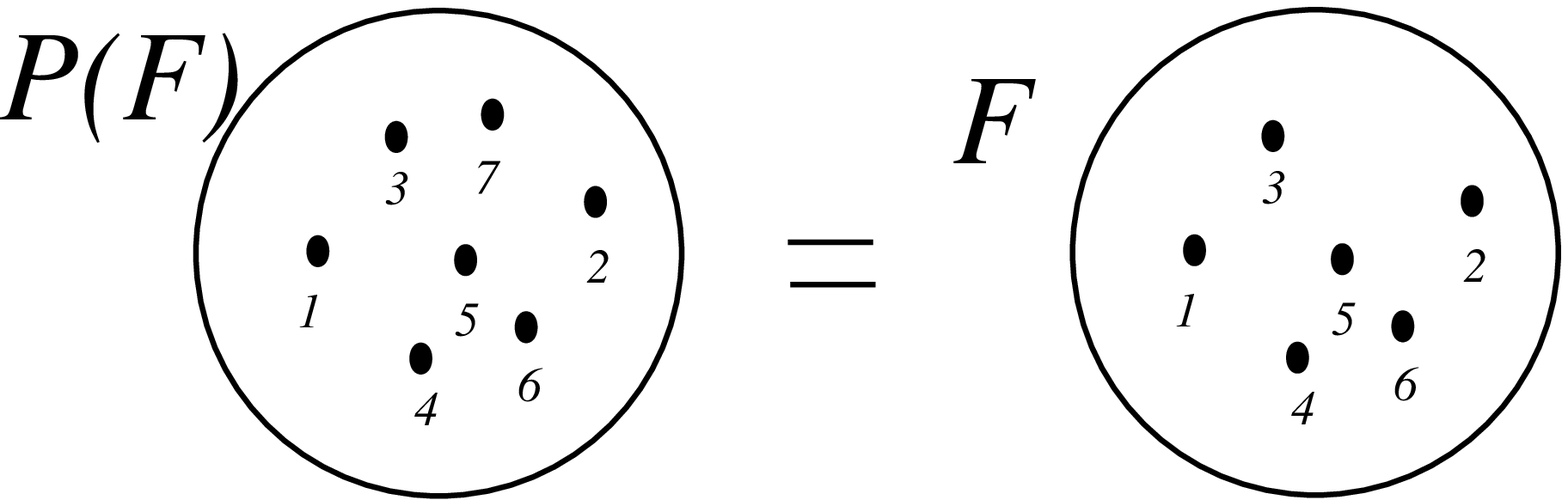}
\caption{ \ Graphical representation of $P$. \label{rb1}}
\end{center}
\end{minipage}
\begin{minipage}[t]{0.5\linewidth}
\begin{center}
\includegraphics[height=2cm]{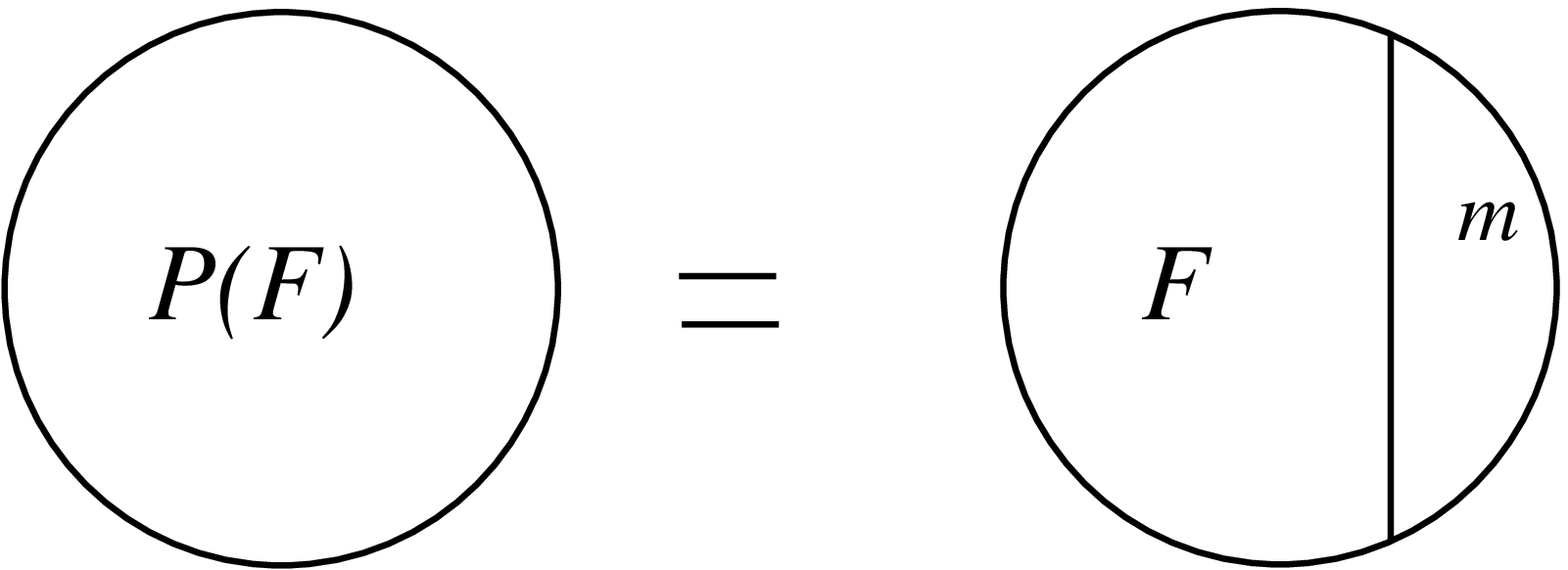}
\caption{ \ Abstract representation. \label{rb2}}
\end{center}
\end{minipage}
\end{figure}

\begin{figure}[ht]
\begin{center}
\includegraphics[height=2cm]{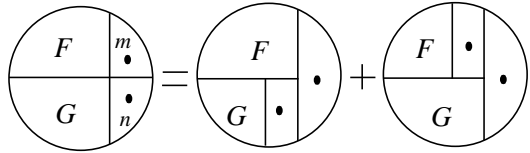}
\caption{ \ Graphical representation of the Rota-Baxter identity. \label{rb3}}
\end{center}
\end{figure}

The action of $P$ on a species $F$ is graphically represented in
Figures \ref{rb1} and \ref{rb2}. With this conventions the
Rota-Baxter isomorphism for $P$ is represented in Figure
\ref{rb3}. Both sides are isomorphic because either $n<m$ and then
the graph on the left agrees with the first graph on the right, or
$m<n$ and in that case it agrees with the second graph on the
right hand side.

\begin{figure}[ht]
\begin{center}
\includegraphics[height=2.5cm]{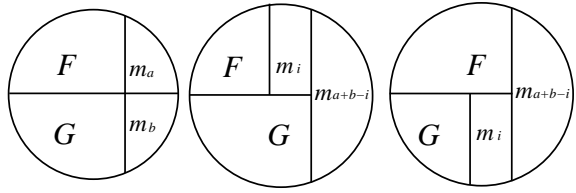}
\caption{ \ Graphical representation of the Rota-Baxter identity. \label{rb4}}
\end{center}
\end{figure}

As in \cite{DPaez} one can show that in any Rota-Baxter category
of weight zero there are natural isomorphisms
$$P^a(x)P^b(y)\simeq \bigoplus_{i=1}^{a}\binom{b-1+a-i}{ b-1 }
P^{a+b-i}(P^i(x)\otimes y) \oplus
\bigoplus_{i=1}^{b}\binom{a-1+b-i}{ a-1 }P^{a+b-i}(F\otimes
P^i(G)).$$ where $a,b > 1$ and $c
\geq 0$ are integers, and by convention $nx$ is the sum of $n$
copies of $x$, for $n$ a non-negative integer and $x$ an object of
a distributive category. It is interesting to elucidate the
meaning of the natural isomorphisms above in the Rota-Baxter
category $(C^{L},P).$ For species $F$ and $G$, the species
$$P^a(F)P^b(G), P^{a+b-i}(P^i(F)\otimes G) \mbox{\ \ and \ \ }
P^{a+b-i}(F\otimes P^i(G))$$ are represented graphically in Figure
\ref{rb4}.  Let us show how the desired isomorphisms arise. Each application of the Rota-Baxter
isomorphisms to $P^a(F)P^b(G)$ yields a couple of the form
$$P^{a+b-1}(P(F)\otimes G) \mbox{\ \  and \ \ }P^{a+b-1}(F\otimes
P(G)).$$ Thus it should be clear that we can apply recursively the
Rota-Baxter isomorphisms until we reach elements of the form
$$P^{a+b-i}(P^i(F)\otimes G) \mbox{\ \ and \ \ } P^{a+b-i}(F\otimes
P^i(G)),$$ where the process stop since then no further
application of the Rota-Baxter isomorphisms is possible. Consider
an application of the functor $P^a(F)P^b(G)$ on a finite set $x$.
The object $P^a(F)P^b(G)(x)$ is sum of a family of objects of $C$
constructed in several steps. First, $x$ is partitioned in two
blocks $x_1$ and $x_2$, then the top $a$, respectively $b$
elements, are remove from $x_1$ and $x_2$, thus we obtain object
$$F(x_1 \setminus m_a(x_1))
\otimes F(x_2
\setminus m_b(x_2)).$$ This case correspond with the left most picture in Figure
\ref{rb4}.  We need to count how many copies of
$$P^{a+b-i}(P^i(F)\otimes G)$$ arise in this process. The species
$P^{a+b-i}(P^i(F)\otimes G)$ applied to a finite set $x$ yields
the sum of the following objects of $C$: first we remove the top
$a+b-i$ elements of $x$, the resulting set is partitioned in two
blocks $x_1'$ and $x_2'$, then we remove the top $i$ elements from
$x_1'$, thus we obtain an object $F(x_1'
\setminus m_i(x_1')) \otimes F(x_2')$. Assume now that the maximal
element of $x$ lies in  $x_2$. The pairs $x_1$ and $x_2$ given
rise to pairs $x_1'$ and $x_2'$ as above are constructed in the
following way: from the $a+b-i$ top elements of $x$ the maximal
element should lie in $x_2$, thus we should choose a subsect $s$
of cardinality $b-1$ from the $a+b-i-1$ top elements (excluding
the maximal element), once this choice has been made $x_1$ and
$x_2$ are uniquely determined from $x_1'$ and $x_2'$ via the
identities $x_1=x_1' \cup m_{a+b-i}(x) \setminus (s \cup
\{m(x)\})$ and $x_2= x_2' \cup s \cup \{m(x)\}.$ Clearly there are
as many as $$\binom{b-1+a-i}{ b-1 }$$ different choices for $s$,
thus justifying the claimed isomorphisms.\\

We closed this section describing a family of Rota-Baxter
categories that may be regarded as categorical analogues of
discrete integration. These examples are based on two techniques
strongly promoted by Rota, the incidence algebra of posets free
Rota-Baxter \cite{GCRota}. The reader will find other approaches
to free Rota-Baxter algebras in
\cite{a}, \cite{cartier}, \cite{guo}. Assume $C$ is a
distributive category and $(X, \leq)$ is a partially order set.
Let $C^{X}$ be the category whose objects are maps $f:X
\rightarrow Ob(C).$ Morphisms in $C^{X}$ from $f$ to $g$ are given by
$$C^{X}(f,g)=\prod_{i \in X} C(f(i),g(i)).$$ Sum and product on
$C^{X}$ are given  by $$(f\otimes g)(i)=f(i)\otimes g(i) \mbox{\ \
and \ \ }(f\otimes g)(i)=f(i)\otimes g(i),$$ respectively. Define
the functors $P_<:C^{X}
\rightarrow C^{X}$  and $P_{\leq}:C^{X}
\rightarrow C^{X}$ by $$P(f)(j)=\oplus_{i < j}f(i) \mbox{\ \ and
\ \ } P(f)(j)=\oplus_{ i \leq j}f(i).$$

\begin{figure}[ht]
\begin{center}
\includegraphics[height=1.5cm]{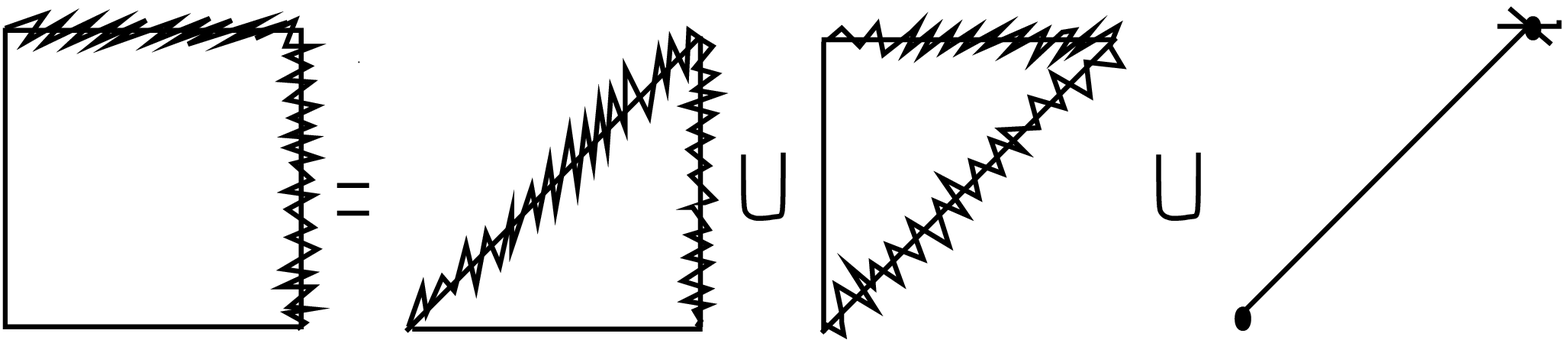}
\caption{ \ Rota-Baxter identity of weight $1$. \label{rb6}}
\end{center}
\end{figure}

\begin{figure}[ht]
\begin{center}
\includegraphics[height=1.5cm]{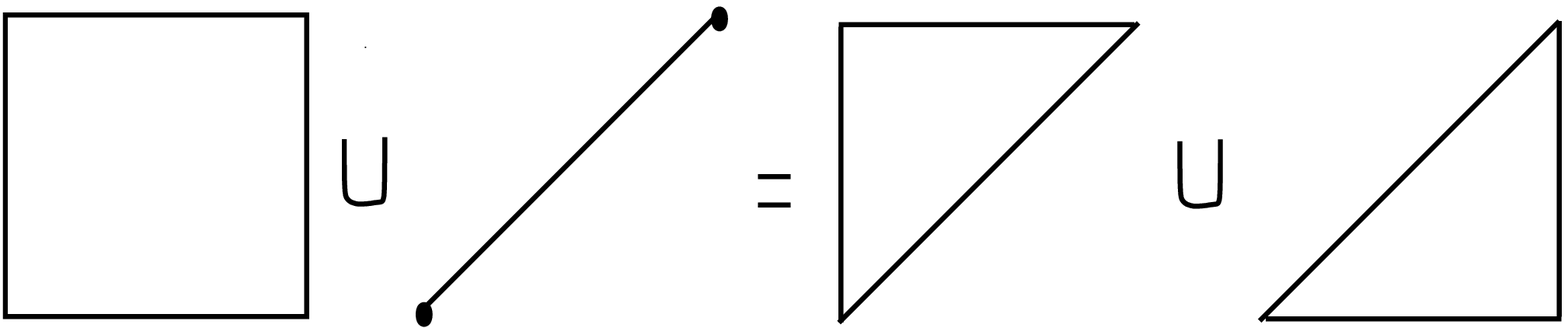}
\caption{ \ Rota-Baxter identity of weight $-1$. \label{rb7}}
\end{center}
\end{figure}

\begin{thm}{\em
$(C^{X},P_<)$ is a Rota-Baxter category of weight $1$.
$(C^{X},P_{\leq})$ is a Rota-Baxter category of weight $-1$.}
\end{thm}
\begin{proof}
The statements follow from the identities between set with
multiplicities depicted in Figure \ref{rb6} and Figure
\ref{rb7}, respectively.
\end{proof}

For any distributive category $C$ and  finite poset $X$ we define
the incidence category $I(X,C)$ as the full subcategory of $C^{X
\times X}$ whose objects $A: X \times X
\longrightarrow C$ are such $A(i,j)=0$ unless $i \leq j$.
$I(X,C)$ is a distributive category with sum and product given by
$(A+B)(i,j)=A(i,j) \oplus B(i,j)$ and $(AB)(i,k)=\bigoplus_{i \leq
j \leq k}A(i,j)\otimes B(j,k)$. There is functor $P:I(X,C)\times
C^X
\longrightarrow C^X$ sending $(A,f)$ into $P_A(f)$ given
by $P_A(f)(j)= \bigoplus_{i \leq j} A(i,j)f(i).$ Letting $\xi$ be
given by $\xi(i,j)=1$ for $i
\leq j$ and $0$ otherwise, then $P_\xi=P_{\leq}$ is a Rota-Baxter operator
of weight $-1$. Similarly, letting $\overline{\xi}$ be given by
$\overline{\xi}(i,j)=1$ for $i < j$ and $0$ otherwise, we obtain
that $P_{\overline{\xi}}=P_{<}$ is a Rota-Baxter operator of
weight $1$. Unfortunately,  $P_A$ is not a Rota-Baxter functor for
arbitrary $A$.

\section{Categorification of $q$-calculus}\label{hoy}

Most applications of $q$-calculus assume that $q$ is a real number
in the interval $(0,1)$. As $q$ approaches $1$ one  recovers
computations in classical calculus. However, it is often the case
that computations in $q$-calculus make sense for $q\geq 0$. In
this section we adopt the convention that $q$ is a non-negative
integer. Setting $q=1$ we recover the theory of species from the
theory of $q$-species developed in this section. Given a
commutative ring $R$ we let $R[[x]]_q$ be the ring of formal
$q$-divided powers series defined as the free $R$-module generated
by symbols $$\frac{x^k}{[k]_q!} \mbox{\ \ with \ \ } k
\in \mathbb{N},$$ with product $$\frac{x^k}{[k]_q!}\frac{x^s}{[s]_q!}=
\binom{k+s}{s}_q \frac{x^{k+s}}{[k+s]_q!},$$ where
$$\binom{k+s}{s}_q=\frac{[k+s]_q!}{[k]_q![s]_q!}.$$ We define
operators $\partial_q: R[[x]]_q \longrightarrow R[[x]]_q$,
$\int_q: R[[x]]_q \longrightarrow R[[x]]_q$ and $s_q: R[[x]]_q
\longrightarrow R[[x]]_q$ by:
$$\partial_q \left( \sum_{n=0}^{\infty}f_n \frac{x^n}{[n]_q!} \right) =\sum_{n=0}^{\infty}f_{n+1} \frac{x^n}{[n]_q!},$$
$$\int_q \left( \sum_{n=0}^{\infty}f_n \frac{x^n}{[n]_q!} \right) =\sum_{n=1}^{\infty}f_{n-1} \frac{x^n}{[n]_q!},$$
$$s_q \left( \sum_{n=0}^{\infty}f_n \frac{x^n}{[n]_q!} \right) =\sum_{n=0}^{\infty}q^nf_{n} \frac{x^n}{[n]_q!}.$$
The operators $\partial_q$ and $\int_q$ are the formal analogues
of the notions of $q$-derivation and $q$-integration. The operator
$s_q$ is called the shift operator and plays a distinguished role
in the our next result. We are going to show that $\int_q$ is a
twisted\footnote{We thank Dominique Manchon for helping us clarify
this subtle point.} Rota-Baxter operator.

\begin{prop} {\em For $f,g \in R[[x]]_q$ the following identity holds:
$$\left(\int_q f\right) \left(\int_q g\right)=  \int_q \left(\int_q f \right) g  +
\int_q  \left( fs_q \left( \int_q g \right) \right) .$$}
\end{prop}
\begin{proof}
Assume that $$f = \sum_{n=0}^{\infty}f_n  \frac{x^n}{[n]_q!}
\mbox{\ \  and  \ \ } g =
\sum_{n=0}^{\infty}g_n  \frac{x^n}{[n]_q!}.$$  The desired result follows
from the identities:
$$(\int_q f)(\int_q g)= \sum_{n=2}^{\infty}\left(\sum_{k=1}^{n-1}\binom{n}{k}_q f_{k-1}g_{n-k-1} \right)\frac{x^n}{[n]_q!},$$

$$ \int_q(\int_q f) g = \sum_{n=2}^{\infty}\left(\sum_{k=1}^{n-1}\binom{n-1}{k}_q f_{k-1}g_{n-k-1} \right)\frac{x^n}{[n]_q!} ,$$

$$\int_q(fs_q(\int_q g) ) = \sum_{n=2}^{\infty}\left(\sum_{k=1}^{n-1}\binom{n-1}{k-1}_q f_{k-1}g_{n-k-1} \right)\frac{x^n}{[n]_q!}   ,$$

$$\binom{n}{k}_q = \binom{n-1}{k}_q + q^{n-k}\binom{n-1}{k-1}_q .$$
\end{proof}

Let $C$ be a distributive category, we are going to define a
$q$-deformed distributive structure on $C^{L_n}$ as follows: sum
of species is given by $(F+G)(x)=F(x)
\oplus G(x)$, and the $q$-product of species is given by $$(FG)(x)=\bigoplus_{x_1
\sqcup x_2 =x}[q]^{c(x_1,x_2)}F(x_1)\otimes G(x_2),$$ where $$c(x_1,x_2)= \{(i,j) \ | \ i \in x_1, j
\in x_2
\mbox{ and } i>j \}.$$  In the definition above we used the following convention: if $x$ is a finite set
and $c$ an object of a distributive category then we set
$$xc=\bigoplus_{i \in x}c.$$   We write $C_q^{L_n}$ instead of
$C^{L_n}$ to emphasize that we are using the $q$-deformed product.
Notice that if $q=1$ then $[q]^{c(x_1,x_2)}$ is a set with one
element and plays no significant role, thus we recover the product
of species of Section \ref{ci}. We define functors $P_q: C_q^{L}
\longrightarrow C_q^{L}$ and $S_q: C_q^{L}
\longrightarrow C_q^{L}$ as follows. For species $F,G$ in $C_q^{L}$ and a linearly order set $x$ we
set $$\partial_q(F)(x)=F(x \sqcup \{x\}),$$
$$P_q(F)(x)=F(x'),$$ where $x'=x \setminus m(x)$ if $x$ is
non-empty and $P_q(F)(\emptyset)=\emptyset$. The functor $S_q$ is
given by $$S_q(F)(x)= [q]^xF(x).$$ One can check that $\partial_q$
and $P_q$ are almost inverse of each other, indeed, we have that
$$\partial_qP_qF=F \mbox{\ \ and \ \ } P_q \partial_q F =F_+,$$ where
$F_+(x)=F(x)$ if $x$ in non-empty and $F_+(\emptyset)= \emptyset.$

\begin{thm}{\em $(C_q^{L}, P_q, S_q)$ is a twisted Rota-Baxter category, i.e. there are natural isomorphisms $$P_q(F)P_q(G)
\simeq P_q(P_q(F)G) \oplus P(FS_qP_q(G)).$$ If $C$ is a
categorification of $R$, then $(C_q^{L}, P_q ,S_q)$ is a
categorification of $(R[[x]]_q,
\int_q, s_q)$.}
\end{thm}

\begin{proof} Let us show that $P_q$ is a twisted Rota-Baxter
functor. We have the identities
$$P_q(F)P_q(G)(x)= \bigoplus_{x_1 \sqcup x_2 =x}[q]^{c(x_1, x_2)} F(x_1') \otimes G(x_2'),$$

$$P_q(P_q(F)G)(x)=\bigoplus_{x_1 \sqcup x_2 =x'} [q]^{c(x_1, x_2)} F(x_1') \otimes G(x_2),$$

$$P_q(FS_qP(G))(x)=\bigoplus_{x_1 \sqcup x_2 =x'}
[q]^{c(x_1, x_2) \sqcup x_2}F(x_1)
\otimes G(x_2'),$$

The desired natural isomorphisms are constructed as follows.
Consider the summand in $P_q(F)P_q(G)(x)$ corresponding with the
partition $x_1 \sqcup x_2=x.$ Then  $m(x)$ lies either in $x_1$ or
in $x_2$. If $m(x) \in x_2$, then $x_1 \sqcup x_2'=x'$ and the
summands corresponding to $(x_1,x_2)$ and $(x_1,x_2')$ in
$P_q(F)P_q(G)(x)$ and $P_q(P_q(F)G)(x)$, respectively, agree since
in this case $c(x_1,x_2)=c(x_1,x_2')$.  On the other hand, if
$m(x) \in x_1$, then $x_1' \sqcup x_2=x'$ and the summands
corresponding to $(x_1,x_2)$ and $(x_1,x_2')$ in $P_q(F)P_q(G)(x)$
and $P_q(FS_qP(G))(x),$ respectively, are naturally isomorphic
since in this case $c(x_1,x_2)=c(x_1',x_2) \sqcup x_2.$\\

The valuation map $| \ \ |: C_q^{L}
\longrightarrow R[[x]]_q$ is given by
$$|F|=\sum_{n=0}^{\infty}|F([n])|\frac{x^n}{[n]_q!}.$$
Let us check that satisfies the multiplicative property.
\begin{eqnarray*}
|FG|&=& \sum_{n=0}^{\infty}|(FG)([n])|\frac{x^n}{[n]_q!} \\
&=& \sum_{n=0}^{\infty}\sum_{x_1
\sqcup x_2 =[n]}q^{|c(x_1,x_2)|}|F(x_1)||G(x_2)|\frac{x^n}{[n]_q!} \\
&=& \sum_{n=0}^{\infty}\sum_{k=0}^{n}\left( \sum_{x_1
\sqcup x_2 =[n], |x_1|=k}q^{|c(x_1,x_2)|} \right) |F(k)||G(n-k)|\frac{x^n}{[n]_q!}\\
&=& \sum_{n=0}^{\infty}\sum_{k=0}^{n}\binom{n}{k}_q |F(k)||G(n-k)|\frac{x^n}{[n]_q!} \\
&=&  |F||G|.
\end{eqnarray*}

\end{proof}

\section{Categorification of classical field theory}\label{cft}

Let $K$ be a set and $J$ a subset of $K$. For a commutative ring
$R$ we let $R[[K]]$ be the ring of formal divided power series on
variables $k \in K$, and $R[[J]]$ the subring of $R[[K]]$
consisting of formal divided power series with variables $k \in
J$. Formally  $R[[K]]$ is defined as follows. Let $M(K)$ be the
set of maps
$$m:K \longrightarrow \mathbb{N}$$ with finite support $s(m)= \{i
\in K \ | \ m(i)\neq 0 \}$, then we set $$R[[K]]= Maps(M(K),R).$$ The
structural operations on $R[[K]]$ are given by $$(f+g)(m)=f(m) +
g(m)$$ and $$(fg)(m)=\sum_{m_1+m_2=m}\binom{m}{m_1}f(m_1)g(m_2),$$
where $m,m_1, m_2
\in M(K)$ and $$\binom{m}{m_1}=\prod_{i \in s(m)}\binom{m(i)}{m_1(i)}.$$ Consider the operator $p_J: R[[K]]\longrightarrow
R[[J]] \subseteq R[[K]]$ given by

$$p_J(f)(m)=\left\{\begin{array}{cc}
f(m) & \mathrm{if}\ s(m) \subseteq x\\
0 & \mathrm{otherwise}
\end{array}\right.$$
It is easy to see that the operator $p_J$ is Rota-Baxter of weight $-1$ since $p_J$ is an idempotent ring morphism.\\

The construction above can be generalized to the categorical
context without difficulties. Let $\mathbb{B}^{(K)}$ be the
category whose objects are pairs $(x,f)$ where $x$ is a finite set
and $f:x \longrightarrow K$ is a map. Morphisms in
$\mathbb{B}^{(K)}$ from $(x,f)$ to $(y,g)$ are bijections
$\alpha:x\rightarrow y$ such that $g \circ
\alpha =f$. Let $C$ be a distributive category with negative objects, see \cite{BD, BD1, BD2, RDEP} for examples; and let
$C^{\mathbb{B}^{(K)}}$ be the category of functors from
$\mathbb{B}^{(K)}$ to $C.$ We define the sum and product of
functors as follows,  given $(x,f)$ in $\mathbb{B}^{(K)}$ and
$F,G$ in $C^{\mathbb{B}^{(K)}}$ then
$$(F+ G)(x,f) = F(x,f)\oplus G(x,f)$$ and $$FG(x,f)=\bigoplus_{x_1
\sqcup x_2 =x}F(x_1,f|_{x_1})\otimes G(x_2,f|_{x_2}).$$ These
structural functors turn $C^{\mathbb{B}^{(K)}}$ into a
distributive category. Assume that $C$ is a categorification of a
ring $R$ and let $R[[K]]$  be the ring of formal series in
variables $K$ with coefficients in $R$. Consider the functor
$P_{J}: C^{\mathbb{B}^{(K)}}
\longrightarrow C^{\mathbb{B}^{(K)}}$ given by
$$P_{J}(F)(x,f)=\left\{\begin{array}{cc} F(x,f) & \mathrm{if}\
f(x)
\subseteq J\\
0 & \mathrm{otherwise}
\end{array}\right.$$
The following result is easy to check.

\begin{thm}\label{foo}{\em
$(C^{\mathbb{B}^{(K)}},P_J)$ is a Rota-Baxter category of weight
$-1$. Moreover $(C^{\mathbb{B}^{(K)}},P_J)$ is a categorification
of $(R[[K]],p_J)$ with valuation map given by $$|F|(m)=|F(x_m,
f_m)|,$$ $$\mbox{\ \ where \ \ }x_m=\coprod_{k \in K}[m(k)]
\mbox{\ \ and  the map \ \ } f_m:\coprod_{k
\in K}[m(k)] \longrightarrow K$$ is such that $f_m(i)= k$ if $i \in [m(k)]$.}
\end{thm}

Let us now see how a Lagrangian field theory may be described from
a categorical point. The basic objects appearing in a field
theory, namely the fields, are locally identified with maps
$$\varphi:
\mathbb{R}^d \longrightarrow \mathbb{R}^k,$$ where $d$ is the dimension of the space-time
manifold and $n$ is the number of scalar fields involved in the
construction of $\varphi.$ A Lagrangian theory is completely
determined by a function $l$, the Lagrangian, depending on fields
and its derivative. Fixing coordinates $x_1,...,x_d$ on
$\mathbb{R}^d$ and writing $\varphi$ as $\varphi =
(\varphi^1,...,\varphi^k),$ then a translation invariant
Lagrangian may be regarded as a polynomial or formal power series
in variables $\partial_I
\varphi^j$, where for $j \in [k]$ and $I \in \mathbb{N}^d$ we set
$$\partial_I
\varphi^j=\partial_1^{I(1)}...\partial_d^{I(d)}\varphi^j.$$
Thus we see that a Lagrangian $l$ may be regarded as an element of
the ring of formal divided powers
$$R[[\mathbb{N}^d \times [k]]]=R[[\partial_I
\varphi^j]],$$ where $(I,j) \in \mathbb{N}^d \times [k]$ and $\partial_I
\varphi^j$ is regarded as a formal variable.
Theorem \ref{foo} tell us that if $C$ is a categorification of
$R$, then $C^{\mathbb{B}^{(\mathbb{N}^d \times [k])}}$ is a
categorification of $R[[\partial_I
\varphi^j]]$. Objects in
$C^{\mathbb{B}^{(\mathbb{N}^d \times [k])}}$ are triples $(x,f,g)$
where $x$ is a finite set, $f: x \longrightarrow \mathbb{N}^d$ and
$g:x \longrightarrow [k]$.  The valuation map $$|
\ \ | :C^{\mathbb{B}^{(\mathbb{N}^d \times [k])}}\longrightarrow
R[[\partial_I
\varphi^j]]$$ sends a species $F \in
C^{\mathbb{B}^{(\mathbb{N}^d \times [k])}}$ into the formal
divided power series
$$|F|=\sum_{f \in M(\mathbb{N}^d \times [k])} F(\bigsqcup_{(I,j)}[(f,g)^{-1}(I,j)]) \frac{\varphi^f}{f!},$$
where $$\frac{\varphi^f}{f!} = \prod_{(I,j)} \frac{(\partial_I
\varphi^j)^{f(I,j)}}{f(I,j)!}.$$ For $(I,j)
\in \mathbb{N}^d \times [k]$ we have the partial derivation functors
$$\partial_{(I,j)}:C^{\mathbb{B}^{(\mathbb{N}^d
\times [k])}} \longrightarrow C^{\mathbb{B}^{(\mathbb{N}^d \times
[k])}}$$  given  by $$\partial_{(I,j)}F(x,f)=F((x,f) \sqcup (*,
(I,j)).$$ Notice that the ring $R[[\partial_I
\varphi^j]]$ comes with additional natural vector fields
$$\partial_i:R[[\partial_I
\varphi^j]] \longrightarrow R[[\partial_I
\varphi^j]]$$ given by
$$\partial_i(\partial_I
\varphi^j)= \partial_{I+\varepsilon_i}\varphi^j$$ where the vectors $\varepsilon_i$ are
the canonical generators of $\mathbb{N}^d$. In the categorical
context we have functors $$\partial_i:C^{\mathbb{B}^{(\mathbb{N}^d
\times [k])}} \longrightarrow C^{\mathbb{B}^{(\mathbb{N}^d \times
[k])}}$$ given for $1 \leq i \leq n$ by $$\partial_i(F)(x,f,g)=
\bigoplus_{a
\in x}F(x, f +
\delta_a\varepsilon_i,g )$$ where $\delta_a: x \longrightarrow \{0,1
\}$ is the Kronecker delta function. More generally for $I \in \mathbb{N}^d$ we have
differential functor $$\partial_{I}:C^{\mathbb{B}^{(\mathbb{N}^d
\times [k])}} \longrightarrow C^{\mathbb{B}^{(\mathbb{N}^d \times [k])}}$$
given by $$\partial_{I}(F)(x,f,g)= \bigoplus_{a_i \in
x^{I(i)}}F(x, f +
\sum_{i=1}^l\delta_{a_{i,j}}\varepsilon_i,g ).$$

The categorification of a Lagrangian $l$ is thus obtained by
finding a functor $L$ in $C^{\mathbb{B}^{(\mathbb{N}^d
\times [k])}}$ such that $|L|=l$.  Classical field theory main concern is understanding
the solutions of a system of partial differential equations
$e_j(l)=0$, $j \in [k],$  called the Euler-Lagrange equations,
determined by the Lagrangian $l$. If we are interested in the
categorification of this system of partial differential equations
the first step is to find categorification for each of the
equations appearing in the Euler-Lagrange equations, namely, we
need species $E_j(l)$ in $C^{\mathbb{B}^{(\mathbb{N}^d
\times [k])}}$ such that $|E_j(L)|=e_j(l).$

\begin{thm}{\em The functors $E_j(L): C^{\mathbb{B}^{(\mathbb{N}^d
\times [k])}} \longrightarrow C^{\mathbb{B}^{(\mathbb{N}^d
\times [k])}}$ given by
$$E_j(L)(x,f,g)= \bigoplus_{I \in \mathbb{N}^d}
\bigoplus_{a_i \in (x\sqcup \{* \})^{I(i)}}(-1)^{|I|}L(x \sqcup \{* \}, (f \sqcup
(*,I))
 + \sum_{i,m}\varepsilon_i \delta_{a_{i,m}}, g \sqcup (*,j))$$ are such that $|E_j(L)|= e_j(L)$.}
\end{thm}

\begin{proof}
Follows from the definitions given above and the well-know formula
$$e_j(l)=\sum_{I \in \mathbb{N}^d}(-1)^{|I|}\partial_I \frac{\partial l}{\partial (\partial_I
\varphi^j)}.$$
\end{proof}

\section{Categorification of quantum field theory}\label{qft}

As discussed in the previous section the ring of functions on the
configuration spaces of a field theory may be identified with
$R[[K]]$ where $K=\mathbb{N}^d \times [k]$. Moreover
$C^{\mathbb{B}^{(K )}}$  provides a categorification of the ring
of functions on configuration space. In order to proceed to
consider quantum structures in field theory, we consider functions
in phase space which may be identified with the ring
$$\mathbb{R}[[K \sqcup
\overline{K}]].$$
An element of the first copy of $K$ is denoted by $k$, the
corresponding element in the second copy is denoted
$\overline{k}$. So we have an involution
$$K
\sqcup \overline{K} \longrightarrow K \sqcup \overline{K}$$ sending $k$ into
$\overline{k}$ and $\overline{k}$ into $k$. Clearly
$C^{\mathbb{B}^{(K
\sqcup \overline{K})}}$ is a categorification of $\mathbb{R}[[K \sqcup
\overline{K}]].$
The new structure algebraic structure present in the ring of
functions on phase space is the Poisson bracket
$$\{\ \ , \
\
\}: \mathbb{R}[[K \sqcup
\overline{K}]] \times \mathbb{R}[[K \sqcup
\overline{K}]]
\longrightarrow \mathbb{R}[[K \sqcup
\overline{K}]]$$
which is completely determined by the fact that it antisymmetric,
a derivation on each variable, and is defined on generators by
$$\{k,{\overline{k}} \}= \delta_{k,\overline{k}}, \ \
\{k, s \}=0, \mbox{ and } \{{\overline{s}} ,{\overline{k}} \}=0.$$

So our first to define a bifunctor
$$\{\ \ , \
\
\}: C^{\mathbb{B}^{(K \sqcup
\overline{K})}} \times C^{\mathbb{B}^{(K \sqcup
\overline{K})}}
\longrightarrow C^{\mathbb{B}^{(K \sqcup \overline{K})}}$$ which plays the role of the Poisson bracket in the categorical
context. The Poisson bracket of functors $F$ and $G$ in
$C^{\mathbb{B}^{(K
\sqcup \overline{K})}}$ turns out to be given by:
\begin{eqnarray*}
\{F,G \}(x,f) &=&\bigoplus_{x_1 \sqcup x_2 =x, k \in K} F(x_1 \sqcup \{* \},f \sqcup (*,k))
\otimes G(x_2 \sqcup \{* \},f \sqcup (*,\overline{k}))\\
&-&\bigoplus_{x_1 \sqcup x_2 =x, k \in K}F(x_1 \sqcup \{*
\},f \sqcup (*,\overline{k}))
\otimes G(x_2 \sqcup \{* \},f \sqcup (*,k)).
\end{eqnarray*}

\begin{thm}{\em For $F, G$ in
$C^{\mathbb{B}^{(K
\sqcup \overline{K})}}$
the following identity holds $$|\{F,G \}|=\{|F|,|G| \}.$$}
\end{thm}

\begin{proof}
The result follows from the fact that the Poisson bracket of
functions on phase space is given by $$\{f,g\}=\sum_{k \in
K}\frac{\partial f}{\partial k}
\frac{\partial g}{\partial {\overline{k}}} -
\frac{\partial f}{\partial {\overline{k}}}
\frac{\partial g}{\partial {k}},$$ and the fact that the formula for the Poisson
bracket of species given above may be also be defined by

$$\{F,G\}=\sum_{k \in K}\partial_k F
\partial_{\overline{k}}G -
\partial_{\overline{k}}F
\partial_k G.$$

\end{proof}
The commutative product of functions on phase space comes equipped with
a natural deformation into a quantum product, often called the Moyal product,
which is determined by the Poisson bracket. Our next goal is to describe the Moyal
product at the categorical level. Recall that quantum phase space possesses an
extra dimension parameterized by a formal variable $\hbar$. Thus a categorification of
formal power series in quantum phase space is given by $C^{\mathbb{B}^{(K
\sqcup
\overline{K} \sqcup \hbar)}}$ with the natural valuation map. Objects in  $\mathbb{B}^{(K
\sqcup \overline{K} \sqcup \hbar)}$ are triples $(x,f,h)$ where $x$ is a finite set, $f:x
\longrightarrow K \sqcup \overline{K}$ is a map, and $h$ is another finite
set. Given a map $\alpha: h \rightarrow K \sqcup \overline{K}$, we
define $s(\alpha)$  the sign of $\alpha$ as  follows: $\pm 1$
according to the parity of $\alpha^{-1}(\overline{K}).$

$$s(\alpha)=\left\{\begin{array}{cc} 1 & \mathrm{if}\
|\alpha^{-1}(\overline{K})| \mathrm{\ \ is \ \  even \ \ }\\
-1 & \mathrm{if}\ |\alpha^{-1}(\overline{K})| \mathrm{\ \ is \
\ odd.
\ \ }
\end{array}\right.$$

We define the Moyal star product of species $F$ and $G$ in
$C^{\mathbb{B}^{(K
\sqcup
\overline{K} \sqcup \hbar)}}$ is  given by
\begin{eqnarray*}
F \star G(x,f,h) &=&\bigoplus s(\alpha)F(x_1
\sqcup h_3,f \sqcup \alpha, h_1)
\otimes G(x_2 \sqcup h_3,f \sqcup \overline{\alpha}, h_2)),
\end{eqnarray*}
where the sum runs over finite sets $x_1, x_2, h_1,h_2, h_3$ such
$x_1 \sqcup x_2 =x$, $h_1 \sqcup h_2 \sqcup h_3 =h$ and $\alpha:
h_3
\rightarrow K \sqcup
\overline{K}.$

\begin{thm}{\em For $F,G$ in $C^{\mathbb{B}^{(K
\sqcup
\overline{K} \sqcup \hbar)}}$ the following identity holds $$|F \star G|=|F| \star |G|.$$}
\end{thm}

\begin{proof}
The result is an instance of a general categorification theorem
for the Kontsevich's star product proved in \cite{RDEP}.
Alternatively, one notices that the expression given for the
$\star$-product of species is the categorical version of the
following expression for the $\star$-product of functions on phase
space $$f\star g =
\sum_{n=0}^{\infty}\frac{h^n}{n!} m( \sum_{k \in K}
\frac{\partial}{\partial k} \otimes \frac{\partial}{ \partial {\overline{k}}} -
\frac{\partial}{ \partial {\overline{k}}} \otimes \frac{\partial}{\partial k})^n(f \otimes g),$$
where for $m$ denotes the product of functions on phase space.
\end{proof}

We close the paper with an example of a quantum-like Rota-Baxter
category. Fix a subset $J$ of $K$. Consider the functor
$P_{J}: C^{\mathbb{B}^{(K
\sqcup
\overline{K} \sqcup \hbar)}}
 \longrightarrow C^{\mathbb{B}^{(K
\sqcup
\overline{K} \sqcup \hbar)}}$ given by

$$P_{J}(F)(x,f,h)=\left\{\begin{array}{cc} F(x,f,h) & \mathrm{if}\
f(x) \subseteq J \sqcup \overline{J}\\
0 & \mathrm{otherwise.}
\end{array}\right.$$

\begin{thm}\label{f}{\em
$(C^{\mathbb{B}^{(K
\sqcup \overline{K} \sqcup \hbar)}},
\star, P_J)$ is Rota-Baxter category of weight $-1$. $(C^{\mathbb{B}^{(K \sqcup \overline{K})}},\star, P_J)$ is a categorification of the Rota-Baxter algebra
$$(R[[K \sqcup \overline{K}]][[ \hbar]], \star, p_J)$$ with
valuation map given by $$|F|(m,n)=|F(x_m, f_m,[n])|,$$ where
$x_m=\coprod_{k
\in K}([m(k)] \sqcup [m(\overline{k})])$ and $f_m:\coprod_{k
\in K}([m(k)] \sqcup [m(\overline{k})]) \longrightarrow K$ is such that

$$f_m(i)=\left\{\begin{array}{cc} k & \mathrm{if}\
i \in [m(k)]\\
\overline{k} & i
\in [m(\overline{k})].
\end{array}\right.$$
}
\end{thm}

\section*{Acknowledgments}
Thanks to Takashi Kimura, Dominique Manchon,  Marcelo P\'aez,
Raymundo Popper, Steven Rosenberg and Sylvie Paycha.

\noindent edmundo.castillo@ciens.ucv.ve\\
Facultad de Ciencias, Universidad  Central de Venezuela, Caracas, Venezuela \\

\noindent ragadiaz@gmail.com\\
Facultad de Administraci\'on, Universidad del Rosario, Bogota, Colombia \\

\end{document}